%% file: rectifible-scs.tex
\documentclass[preprint]{elsarticle}
\input in/rtg-ru-head.tex

\begin{document}

\begin{frontmatter}

\title{Algebraic structures on the Cantor set}
\author{Evgenii Reznichenko}

\ead{erezn@inbox.ru}

\address{Department of General Topology and Geometry, Mechanics and  Mathematics Faculty, 
M.~V.~Lomonosov Moscow State University, Leninskie Gory 1, Moscow, 199991 Russia}

\begin{abstract}
Below, by space we mean a separable metrizable zero-dimensional space.
It is studied when the space can be embedded in a Cantor set while maintaining the algebraic structure.
Main results of the work: every space is an open retract of a Boolean precompact group; every strongly homogeneous space is rectifiable. In this case, the space can be embedded in the Cantor set with the preservation of the algebraic structure.
An example of a strongly homogeneous space is constructed which do not admit the structure of a right topological group.
\end{abstract}
\end{frontmatter}

\def\hm{\sim}
\def\hm{\approx}
\def\pp{\Phi}
\def\fuv#1{\Theta\left[#1\right]}
\def\ddt{\D^{<\om}}
\def\cC{\mathcal{C}}
\def\cCs{\cC^*}
\def\smz/{SMZD}
\def\lla{\Lambda}
\def\llan{\lla_*}
\def\opx{\Gamma}
\def\rr{\Upsilon}

\def\P{{\mathbf{P}}}
\def\Q{{\mathbf{Q}}}

\def\T{{\mathbf{T}}}
\def\S{{\mathbf{S}}}

\input sec/intro.tex

\input sec/defs.tex

\input sec/nas.tex

\input sec/rexscs.tex

\input sec/rexscsc.tex

\input sec/hsnrtg.tex

\input sec/qe.tex

\bibliographystyle{elsarticle-num}
\bibliography{rectifible-scs}
\end{document}

%% file: in/rtg-ru-head.tex
\usepackage{cmap} 
           \usepackage[T2A]{fontenc}   
           \usepackage[utf8]{inputenc}
           \usepackage[english]{babel}
\usepackage{amsmath}
\usepackage{amsfonts}
\usepackage{amssymb}
\usepackage{amsthm}
\usepackage[pdftex,unicode]{hyperref}
\usepackage{amscd}
\usepackage[all,cmtip]{xy}
\usepackage{mathtools}
\usepackage{comment}
\usepackage{tikz-cd}
\usetikzlibrary{babel}

\def\T{\mathbb T}  
\def\C{{\mathbf{C}}}
\def\D{{\mathbf{D}}}

\def\om{\omega}

\def\aut#1{\operatorname{Aut}(#1)}
\def\supp{\operatorname{supp}}

\def\be{\beta}

\def\al{\alpha}
\def\ga{\gamma}

\def\nfs/{NFS}
\def\cdp/{CDP}
\def\cdpz/{CDP${}_0$}

\def\cG{\mathcal{G}}

\def\sq#1#2{(#1)_{#2}}
\def\sqn#1{\sq{#1}{n\in\om}}
\def\sqnn#1{\sqn{#1_n}}

\def\cl#1{\overline{#1}}

\def\id{\operatorname{id}}

\def\es{\varnothing}

\def\nom{{n\in\om}}

\def\sset#1{\{#1\}}

\def\set#1{\bbset#1\eeset}
\def\bbset#1:#2\eeset{\{#1\,:\,#2\}}

\def\bbsett#1:#2\eesett{\{#1\,:\,\text{#2}\}}

\def\ibbset#1:#2\ieeset{(#1)_{#2}}

\def\cP{{\mathcal P}}
\def\cB{{\mathcal B}}

\def\cB{{\mathcal B}}

\newcommand\restrA[2]{{
  \left.\kern-\nulldelimiterspace 
  #1 
  \vphantom{\big|} 
  \right|_{#2} 
  }}

\newcommand\restrB[2]{\ensuremath{\left.#1\right|_{#2}}}

\def\restr#1#2{\restrB{#1}{#2}}

\def\pwr#1_#2{#1^{[#2]}}
\def\term#1{{\it #1}}

\def\pp{{\mathsf{P}}}

\def\dddm#1(#2){N_{#1}(#2)}
\def\dddb#1(#2){B_{#1}(#2)}

\def\wh#1{\widehat{#1}}

\def\et(#1){ (#1)}

\def\bitem#1,#2.{ $#2\nrightarrow #1$:\ }
\long\def\edemo#1\endedemo{{\rm #1}}

\newtheorem{assertion}{Statement}

\newtheorem{proposition}{Proposition}
\newtheorem{theorem}{Theorem}
\newtheorem*{theorem*}{Theorem}
\newtheorem{lemma}{Lemma}
\newtheorem*{lemma*}{Lemma}
\newtheorem{cor}{Corollary}

\theoremstyle{definition}
\newtheorem{example}{Example}
\newtheorem{definition}{Definition}
\newtheorem{question}{Question}

\theoremstyle{remark}
\newtheorem{note}{Remark}
\newtheorem*{note*}{Remark}

\def\oo#1/{$O_{#1}$}

\def\rr{$R$}



%% file: sec/intro.tex
\section{Introduction}

In this paper, we study what algebraic structures are possible on subspaces of the Cantor set $\C$, that is, on separable metrizable zero-dimensional (\smz/) spaces. Any \smz/ space is Mal'tsev, moreover, it is a retract of the topological group, that is, a retral space \cite{grs1997}.

We are primarily interested in Mal'tsev, rectifiable, retral, homogeneous \smz/ spaces, (right) topological groups.
It is also studied when a \smz/  space with algebraic structure can be embedded in \smz/ a compact space with the same structure. In particular, when the \smz/ space can be embedded in $\C$ with the algebraic structure preserved.

Main results of the work: each \smz/ space is an open retract of a Boolean precompact \smz/ group (Theorem \ref{tnas1}); every \smz/ strongly homogeneous space is rectifiable (Theorem \ref{trexscsc-main}).

Rectifiable spaces are also characterized as spaces on which there is a Mal'tsev operation with some additional identity ($\pp_3$), the so-called homogeneous Mal'tsev operations (Definition \ref{drexscs1}).
A strongly Mal'tsev operation is introduced and studied, which defines strongly rectifiable spaces. This class lies between rectifiable spaces and topological groups.

In the \ref{sec-hsnrtg} section, an example of a strongly homogeneous \smz/ space is constructed which do not admit the structure of a right topological group.

%% file: sec/defs.tex
\section{Definitions and notation}

By a space we understand a regular topological space.

The weight and character of the space $X$ will be denoted by $w(X)$ and $\chi(X)$, respectively.

Denote by $\aut X$ the set of all homeomorphisms of the space $X$ onto itself, $\aut{X|M}=\set{f\in\aut X: f(M)=M}$ for $M\subset X $. Denote by $\id_X$ the identity mapping of $X$ onto itself.

A space $X$ is called \term{strongly homogeneous} if every non-empty clopen subset of $X$ is homeomorphic to $X$. Every strongly homogeneous \smz/ space is homogeneous.

If $f:X\to Y$ is a mapping and $A\subset X$, then $\restr fA$ denotes the restriction of $f$ to $A$.

\subsection{Non-Archimedean spaces}

A subset $M$ of a linearly ordered set $X$ is called convex if $a\leq b\leq c$ holds for $a,c\in M$ and $b\in X$, then $b\in M$.

A family of sets $\cB$ is called \term{non-Archimedean} if for $U,V\in\cB$ the following holds: if $U\cap V\neq\es$, then either $U\subset V$ or $ V\subset U$. \term{Non-Archimedean base} is the base $\cB$ of the space $X$, which is a non-Archimedean family. Non-Archimedean bases called bases of rank 1 were introduced and studied in \cite{arhangelskii1962,arhangelskii1963}.

Every non-Archimedean space has a base which is a tree by reverse inclusion (Theorem 2.9 \cite{Nyikos1999}).
Each a base which is a tree by reverse inclusion is a non-Archimedean base.

A linear order $\leq$ on a space $X$ is said to be \term{consistent with a non-Archimedean base $\cB$} if any element of $U\in\cB$ is convex with respect to the order $\leq$.

\subsection{Cantor set}

Denote $\D=\sset{0,1}$, $\ddt=\bigcup_\nom \D^n$. Here $\D^0=\sset\es$.

Denote by $\C$ the Cantor set, $\C=\D^\om$.
Let $x=\sqnn x\in\C$, $k\in\om$ and $c=\sq{c_n}{n<k}\in\D^k \subset \ddt$. Denote
\begin{align*}
x|k&=(x_0,...,x_{k-1}), &
x(k)&=x_k, &
B(c)&=\set{x\in\C:x|k=c}.
\end{align*}
The family \[\cC=\set{B(c):c\in \ddt}\] is the standard non-Archimedean base $\C$.
The base $\cC$ is a tree by reverse inclusion.
 On $\C$ we will consider the lexigraphic order: for different $x,y\in \C$, $x<y$ if $x|k=y|k$, $x(k)=0$ and $y( k)=1$ for some $k\in\om$. The lexigraphic order is consistent with the non-Archimedean base $\cC$.

Denote $\supp x=\set{\nom: x(n)=1}$ for $x\in\C$, $\Q=\set{x\in\C: |\supp x|<\om }$ and $\P=\C\setminus \Q$. The space $\Q$ is homeomorphic to rational numbers, and the space $\P$ is homeomorphic to irrational numbers. 

\subsection{Rectifiable, Maltsev, retral spaces}

Let $X$ be a space, $e\in X$, $\Psi: X^2\to X^2$. A mapping is called \term{rectiﬁcation} if it is a homeomorphism and there exists a mapping $p: X^2\to X$ such that $\Psi(x,y)=(x,p(x,y))$ and $p (x,x)=e$ for $x,y\in X$.
A space $X$ is called \term{rectiﬁable} if there exists a rectiﬁcation on $X$.
We define the operation $q: X^2\to X$ by the identity $q(x,p(x,y))=y$. The operation $q$ is continuous and the identities $q(x,e)=x$ and $p(x,q(x,y))=y$ hold for it. Also, $\Psi^{-1}(x,y)=(x,q(x,y))$.
The operations $p$ and $q$ are called \term{homogeneous algebra} \cite{choban1992,ArhangelskiiChoban2010}.

A biternary algebra on a space $X$ is a pair of ternary continuous operations $\al,\be: X^3 \to X$ such that
\[
 \al (y,y,x)=\al(\be(x,y,z),y,z)=\be(\al(x,y,z),y,z)=x,
\]
for all $x,y,z \in X$ (see \cite{malcev1954}).

For a space $X$ the following conditions are equivalent: $X$ is a rectiﬁable space;
$X$ is homeomorphic to a homogeneous algebra; there exists a structure of a biternary algebra on $X$ \cite{choban1992}.

A mapping $\pp: X^3\to X$ is called a \term{Maltsev operation} if $\pp(x,y,y)=\pp(y,y,x)=x$ for $x,y\in X$.
A space $X$ is called \term{Mal'tsev} if there exists a continuous Mal'tsev operation on $X$.
Rectifiable spaces are Maltsev: $\pp(x,y,z)=q(x,p(y,z))$.

A space $X$ is called \term{retral} if there exists a topological group $G$ into which $X$ is embedded in such a way that $X$ is a retract of $G$ \cite{arkhangelskii1987}. The retral space is Mal'tsev, the reverse is not true \cite{grs1997}.

%% file: sec/nas.tex
\section{Non-Archimedean spaces}\label{sec-nas}

In this section, we prove a more detailed and specialized version of Theorem 7 (version 2) from \cite{grs1997} and apply it to $\C$.

Let $X$ be a zero-dimensional space. Denote by $B(X)$ the set of all finite subsets of $X$.
We consider $B(X)$ as an adiative Boolean group with zero $\es$ and the group operation of addition is a symmetric difference.
Denote by $\opx(X)$ the family of all open partitions of $X$. For $\ga\in \opx(X)$ we set
\[
H(\ga)=\set{g\in B(X): |g\cap U|\text{ is even for }U\in\ga}.
\]
The set $H(\ga)$ is a subgroup of $B(X)$. Let $\cG\subset \opx(X)$. Denote by $B_z(X,\cG)$ the group $B(X)$ with the group topology whose prebase at unity is formed by sets of the form $H(\ga)$ for $\ga\in\cG$. The set $X$ naturally embeds in $B(X)$: $i: X\to B(X),\ x\mapsto \sset x$. This embedding is topological and $i(X)$ is closed in $B(X)$ if and only if $\bigcup \cG$ is a subbase of $X$. In what follows we will identify $X$ and $i(X)$.

Denote by $B_z(X)=B_z(X,\opx(X))$. The space $X$ is close embedded in $B_z(X)$. Let $Y\subset X$. Denote
\begin{align*}
B_z^0(X)&=\set{g\in B_z(X): |g|\text{ is even}},
&
B_z^1(X)&=\set{g\in B_z(X): |g|\text{ is odd}},
\\
B_z(X|Y)&=\set{g\in B_z(X): g\subset Y},
\\
B_z^0(X|Y)&=B_z(X|Y) \cap B_z^0(X),
&
B_z^1(X|Y)&=B_z(X|Y) \cap B_z^1(X).
\end{align*}
The set $B_z^0(X)=H(\sset X)$ is a clopen subgroup of index two.
The set $B_z^1(X)=B_z(X)\setminus B_z^0(X)$ is openly closed in $B_z(X)$ and is a coset of the group $B_z^0(X)$ \cite{grs1997}.

Let $X$ be a space with non-Archimedean base $\cB$ is a tree by reverse inclusion and some linear order $<$ consistent with base $\cB$.

For $\cB'\subset \cB$ denote by $\llan(\cB')$ the set of maximal by inclusion elements $\cB'$. Then
$\llan(\cB') \subset \cB'$;
the family $\cB'$ is inscribed in $\llan(\cB')$;
$\llan(\cB')$ is an open partition of $\bigcup \cB'$.

For a family $\ga$ of open subsets $X$, denote $\lla(\ga)=\llan(\cB')$, where $\cB'=\set{U\in\cB: U\subset V\text { for some }V\in\ga}$. Then $\lla(\ga)$ is inscribed in $\ga$; $\lla(\ga)$ is an open partition of $\bigcup\ga$; if $V\in\cB$ and $V\subset W$ for some $W\in\ga$, then $V\subset U$ for some $U\in \lla(\ga)$. If $\ga$ is an open cover of $X$, then $\lla(\ga)\in\opx(X)$ and $\lla(\ga)$ are inscribed in $\ga$. Therefore, $B_z(X)=B_z(X,\cP)$, where
\[
\cP=\set{\ga\in \opx(X): \ga \subset \cB}.
\]

For an open $W\subset X$, denote $\lla(W)=\lla(\sset W)$. Then $\lla(W)$ is an open partition of $W$; if $V\in\cB$ and $V\subset W$, then $V\subset U$ for some $U\in \lla(W)$.

Recall the construction of the retraction $r: B_z^1(X)\to X$ from \cite{grs1997}. For $g\in B_z^1(X)$ we set $\rr(g)=\llan(\ga)$, where $\ga=\set{U\in\cB: g\cap U\neq\es \text{ and }|g\cap U|\text{ is even}}$. Let $R(g)=g\setminus \bigcup \rr(g)$ and $r(g)=\min R(g)$. Theorem 7 (version 2) \cite{grs1997} proves that $r$ is a continuous retraction.

\begin{lemma}\label{lnas1}
The mapping $r$ is open. If $Y\subset X$ and $S=B_z^1(X|Y)$, then the mapping
$\restr rS: S\to Y$ is a continuous open retraction of $S$ to $Y$.
\end{lemma}
\begin{proof}
Let $g\in B_z^1(X)$ and $W\subset B_z(X)$ be a neighborhood of the point $g$. Then $g+H(\ga)\subset W$ for some $\ga\in \opx(X)$. Let $x=r(g)$. Let $x\in V \in \ga$. The family $\rr(g)$ is finite and $x\notin \bigcup \rr(g)$. Take $U\in\cB$ such that $x\in U \subset V\setminus \bigcup \rr(g)$.

Let $y\in U$ and $h=g + \sset{x,y}$. Then $\rr(g)=\rr(h)$, $R(h)=R(g) + \sset{x,y}$ and $r(h)=y$.
Let $M=\set{g + \sset{x,y}: y\in U}$. Then $M\subset g+H(\ga)$ and $r(M)=U$. Hence $r(g)\in U \subset r(W)$.

Let $g\in Y$, $Q=U\cap Y$ and $L=\set{g + \sset{x,y}: y\in Q}$. Then $L=M\cap S$ and $r(L)=Q$. Hence $r(g)\in Q \subset r(W\cap S)$.
\end{proof}

Lemma \ref{lnas1} and Theorem 7 (version 2) \cite{grs1997} imply the following assertion.

\begin{theorem}\label{tnas1}
Let $X$ be a non-Archimedean space. Then there exists a continuous open retraction $r: B_z^1(X)\to X$ such that for any $Y\subset X$ and $S=B_z^1(X|Y)$ the map $\restr rS:S\to Y$ is a continuous open retraction.
\end{theorem}

If $X$ is a compact space, then $B_z(X)$ is a precompact group. The set $\opx(X)$ is countable. Hence $B_z(\C)$ is a precompact \smz/ group.
Let's apply the theorem \ref{tnas1} to $\C$.

\begin{theorem}\label{tnas2}
There exists a continuous open retraction $r: B_z^1(\C)\to \C$ such that for any $X\subset \C$ and $S=B_z^1(\C|X)$ the map $\restr rS: S\to X$ is a continuous open retraction.

The set $B_z^1(\C)$ is a clopen subset of the precompact \smz/ group $B_z(\C)$. Let $y\in Y$ and
\[
h_y: B_z^1(\C)\to B_z^0(\C),\ g \mapsto g + y.
\]
The mapping $h_y$ is a homeomorphism and $h_y(B_z^1(\C|X))=B_z^0(\C|X)$, i.e. $X$ is an open continuous retract of the precompact Boolean \smz/ group $B_z^ 0(\C|X)$.
\end{theorem}

\begin{cor}\label{cnas1}
Any \smz/ space is a continuous open retract of a precompact Boolean group.
\end{cor}

%% file: sec/rexscs.tex
\section{Mal'tsev and rectifiable spaces}\label{sec-rexscs}

Let $X$ be a space with the operation $\pp:X^3\to X$.
Consider the following conditions for the operation $\pp$:
\begin{itemize}
\item[{\rm ($\pp_1$)}]
$\pp(x,y,y)=\pp(y,y,x)=x$;
\item[{\rm ($\pp_2$)}]
$\pp(x,y,\pp(y,z,u))=\pp(x,z,u)$;
\item[{\rm ($\pp_3$)}]
$\pp(x,y,\pp(y,x,u))=u$.
\end{itemize}
for $x,y,z,u\in X$.
The $\pp$ operation is a Mal'tsev operation if ($\pp_1$) is satisfied.

\begin{definition}\label{drexscs1}
Let $X$ be a space and $\pp:X^3\to X$ be a continuous map.
We call a mapping $\pp$ \term{homogeneous Mal'tsev operation} if $\pp$ satisfies the conditions {\rm ($\pp_1$)} and {\rm ($\pp_3$)}.
We call a mapping $\pp$ a \term{strong Mal'tsev operation} if the conditions {\rm ($\pp_1$)} and {\rm ($\pp_2$)} are satisfied for $\pp$. A space $X$ is said to be \term{strongly rectifiable} if there exists a continuous strong Mal'tsev operation on $X$.
\end{definition}

\begin{proposition}\label{prexscs1}
The space $X$ is rectifiable if and only if there exists a homogeneous Mal'tsev operation $\pp$ on $X$.
\end{proposition}
\begin{proof}
Assume that $X$ is a rectifiable space. Let the operations $p$ and $q$ define the structure of a homogeneous algebra, $\Psi$ is a rectification, $\Psi(x,y)=(x,p(x,y))$. We put $\pp(x,y,z)=q(x,p(y,z))$. $q(x,e)=x$ and $p(x,x)=e$ imply $\pp(x,y,y)=x$. $q(x,p(x,y))=y$ implies $\pp(y,y,x)=x$. Let's check ($\pp_3$).
Since $p(y,q(y,p(x,u)))=p(x,u)$ and $q(x,p(x,u))=u$, then
\[
\pp(x,y,\pp(y,x,u))=q(x,p(y,q(y,p(x,u))))=u.
\]

Suppose that $\pp$ satisfies the conditions ($\pp_1$) and ($\pp_3$). Take $e\in X$ arbitrarily.
Let $p(x,y)=\pp(e,x,y)$ and $q(x,y)=\pp(x,e,y)$. From ($\pp_1$) follows $p(x,x)=e$ and $q(x,e)=x$.
From ($\pp_3$) it follows
\begin{align*}
q(x,p(x,y))&=\pp(x,e,\pp(e,x,y))=y,
\\
p(x,q(x,y))&=\pp(e,x,\pp(x,e,y))=y.
\end{align*}
\end{proof}

\begin{proposition}\label{prexscs2}
If $X$ is strongly rectifiable, then $X$ is rectifiable.
\end{proposition}
\begin{proof}
The condition ($\pp_3$) follows from ($\pp_1$) and ($\pp_2$), it suffices to substitute $z=x$ into ($\pp_3$). It remains to apply the Proposition \ref{prexscs1}.
\end{proof}

\begin{proposition}\label{prexscs2+1}
Let $Y$ be a Hausdorff space, $X\subset Y=\cl X$, mapping $\wh\pp: Y^3\to Y$ is continuous, and $\pp(X^3)=X$.
Denote $\pp=\restr{\wh\pp}{X^3}$.
\begin{enumerate}
\item
If $\pp$ is a Mal'tsev operation, then $\wh\pp$ is a Mal'tsev operation.
\item
If $\pp$ is a homogeneous Mal'tsev operation, then $\wh\pp$ is a homogeneous Mal'tsev operation.
\item
If $\pp$ is a strong Mal'tsev operation, then $\wh\pp$ is a strong Mal'tsev operation.
\end{enumerate}
\end{proposition}
\begin{proof}
If the identities ($\pp_1$), ($\pp_2$), and ($\pp_3$) hold on a dense subspace, then they hold everywhere.
\end{proof}

\begin{proposition}\label{prexscs2+2}
Let $X$ be a space with continuous Mal'tsev operation $\pp$, $Y_1$ and $Y_2$ be compact Hausdorff extensions of $X$, the operation $\pp$ extends to $Y_1$ and $Y_2$ up to continuous operations $\pp_1$ and $\pp_2$, respectively. Then $Y_1$ and $Y_2$ coincide, that is, there exists a homeomorphism $f:Y_1\to Y_2$ for which $\restr f{X}=\id_X$.
\end{proposition}
\begin{proof}
We identify $X$ with $\set{(x,x)\in Y_1\times Y_2: x\in X}$. Let $Y_3$ be the closure of $X$ in $Y_1\times Y_2$, $\pi_i$ be the projection of $Y_1\times Y_2$ onto $Y_i$, $f_i=\restr{\pi_i}{Y_3}$ for $i =1.2$, $\pp_3=\restr{\pp_1\times\pp_2}{X^3}$. Then $\pp_3(X^3)=X$ and, by the proposition \ref{prexscs2+1}, the operations $\pp_i$ for $i=1,2,3$ are a Mal'tsev operation and $f_i$ are morphisms with respect to them. Since the spaces $Y_i$ are compact, the mappings $f_i$ are factorial. Therefore, mappings $f_i$ are open (Mal'tsev's
Theorem 4.11 \cite{ReznichenkoUspenskij1998}, \cite{malcev1954}). An open continuous map that is a homeomorphism on a dense set is itself a homeomorphism. Hence the mappings $f_i$ are homeomorphisms. Therefore, the mapping $f=f_1\circ f_2^{-1}: Y_1\to Y_2$ is a homeomorphism.
\end{proof}

\begin{note}\label{nrexscs1}
Propositions \ref{prexscs2+1}, \ref{prexscs2+2} and their proofs remain valid if we consider separately continuous Mal'tsev operations.
\end{note}

\begin{definition}\label{drexscs2}
Let $X$ be a space and $\pp:X^3\to X$ be a continuous Mal'tsev operation. An operation $\pp$ is called \term{precompact} if there exists a Hausdoff compact extension $bX$ of the space $X$ onto which the mapping $\pp$ extends to the mapping $\wh \pp: bX^3\to bX$.

A compact extension $bX$ is uniquely defined up to isomorphism (Proposition \ref{prexscs2+2}), we denote it by $\wh X$ or $\wh{(X,\pp)}$.
\begin{enumerate}
\item
A space $X$ is called \term{precompact Mal'tsev} if there exists a continuous precompact Mal'tsev operation on $X$.
\item
A space $X$ is called \term{precompact rectifiable} if there exists a continuous precompact homogeneous Mal'tsev operation on $X$.
\item
We call a space $X$ \term{precompact strongly rectifiable} if there exists a continuous precompact strong Mal'tsev operation on $X$.
\end{enumerate}
\end{definition}

\begin{proposition}\label{prexscs3}
If $X$ is a topological group, then $X$ is strongly rectifiable.

If $X$ is a precompact topological group, then $X$ is precompact and strongly rectifiable.
\end{proposition}
\begin{proof}
We put $\pp(x,y,z)=xy^{-1}z$. The operation $\pp$ is the standard Mal'tsev operation on a group. Let's check ($\pp_2$).
\[
\pp(x,y,\pp(y,z,u))=xy^{-1}(yz^{-1}u)=xz^{-1}u=\pp(x,z,u ).
\]

If $G$ is precompact, then $\pp$ extends to the completion of $G$.
\end{proof}

\begin{proposition}\label{prexscs4}
Let $X$ be a Mal'tsev precompact space. Then
$w(\wh X)=w(X)=\chi(X)$.
\end{proposition}
\begin{proof}
Let $\tau=\chi(X)$.
The compact Mal'tsev space $\wh X$ is a Dugundzhi compact (Theorem 1 \cite{uspenskij1989}).
The Dugunji compactum is a dyadic compactum.
If a compact dyadic compact contains a dense subspace whose character does not exceed $\tau$, then its weight does not exceed $\tau$ (\cite{efimov1965}).
\end{proof}

\begin{lemma}\label{lrexscs1}
Let $\pp$ be a precompact homogeneous Mal'tsev operation on the space $X$ and $e\in X$.
For any non-empty open $U\subset X$ there exists a finite set $M\subset X$ so that
$X=\bigcup_{y\in M} \pp(e,y,U)$.
\end{lemma}
\begin{proof}
Let $W\subset \wh X$ be open and $W\cap X=U$. Take a non-empty open $V\subset \cl V \subset W$. The family $\set{\pp(e,y,V):y\in\wh X}$ forms an open cover of $\wh X$. Then $\bigcup_{i=1}^n \pp(e,x_i,V)$ for some $x_1,x_2,...,x_n\in \wh X$. Then $\pp(e,x_i,V)\subset \pp(e,y_i,U)$ for some $y_i\in X$. Let $M=\sset{y_1,y_2,...,y_n}$.
\end{proof}

%% file: sec/rexscsc.tex
\section{Mal'tsev operations on the Cantor set $\C$}\label{sec-rexscsc}

\begin{theorem}\label{trexscsc1}
Let $X$ be the \smz/ space. Then $X$ is a pre-compact Mal'tsev space.

Moreover, there is a precompact Mal'tsev operation $\pp$ on $\C$, so that $\pp(X^3)=X$ for any $X\subset \C$.
\end{theorem}
\begin{proof}
Let $r: B_z^1(\C)\to \C$ as in Theorem \ref{tnas2}. Let $\pp(x,y,z)=r(\sset{x,y,z})$. Then $\pp(x,y,z)\in\sset{x,y,z}$ and hence $\pp$ is a Mal'tsev operation.
\end{proof}

\begin{theorem}\label{trexscsc-main}
Let $X$ be a strongly homogeneous \smz/ space. Then $X$ is a precompact strongly rectifiable space.

Moreover, there exists a precompact strongly Malcev operation $\pp$ on $X$, so that $\wh{(X,\pp)}$ is homeomorphic to $\C$.
\end{theorem}
\begin{proof}
There is a compact zero-dimensional extension $bX$ of the space $X$ homeomorphic to $\C$, so that for any non-empty open-closed $U,V\subset bX$, $U,V\neq Y$ there exists $f\in\aut{ Y|X}$ so $f(U)=V$, Lemma 2.2 \cite{vanMill2007}.
We will assume that $bX=\C$.

\begin{lemma}\label{lrexscsc-1}
For every non-empty open-closed $U\subset \C$ there exists a homeorphism $f_U:\C\to U$, so that $f_U(X)=U\cap X$.
\end{lemma}
\begin{proof}
Let $V_1$ and $V_2$ be a partition of $\C$ into two non-empty open-closed subsets and $U_1$ and $U_2$ be a partition of $U$ into two non-empty open-closed subsets. For $i=1,2$, let $f_i: V_i\to U_i$ be a homeorphism such that $f_i(V_i\cap X)=U_i\cap X$. Then $f_U$ is the union of $f_1$ and $f_2$.
\end{proof}
Denote $\cCs=\cC\setminus \sset\C$.
For $U,V\in\cCs$ we set
\[
\fuv{U,V}=f_U\circ f_V^{-1}: V\to U.
\]
It follows from the construction
\begin{itemize}
\item[{\rm ($\Theta_1$)}]
$\fuv{U,V}: V\to U$ is a homeomorphism and $\fuv{U,V}(V\cap X)=U\cap X$ for $U,V\in\cCs$;
\item[{\rm ($\Theta_2$)}]
$\fuv{U,U}=\id_U$ for $U\in\cCs$;
\item[{\rm ($\Theta_3$)}]
$\fuv{U,V} \circ \fuv{V,W} =\fuv{U,W}$ for $U,V,W\in\cCs$.
\end{itemize}

Let $x=\sqnn x\in\C$, $k\in\om$.
Let's put
\begin{align*}
U_k(x)&=\set{y\in\C:x|k+1=y|k+1}=B(x|k+1),
\\
V_k(x)&=\set{y\in\C:x|k=y|k\text{ and }x(k)\neq y(k)}=U_{k-1}(x)\setminus U_k(x)=B(\tilde x),
\end{align*}
where $\tilde x=(x_0,...,x_{k-1},1-x_k)$. Note that
\begin{itemize}
\item[{\rm ($U_1$)}]
$U_k(x),V_k(x)\in\cCs$;
\item[{\rm ($U_2$)}]
the family $\set{V_n(x): n\leq k}$ is an open partition of the set $\C\setminus U_k(x)$;
\item[{\rm ($U_3$)}]
the family $\set{V_n(x): \nom }$ is an open partition of the set $\C\setminus \sset x$;
\item[{\rm ($U_4$)}]
the family $\set{U_n(x): \nom }$ is the base at the point $x$;
\item[{\rm ($U_5$)}]
$U_k(x)=U_k(y)$, $V_k(x)=V_k(y)$ and $x\in U_k(y)$ if $y\in U_k(x)$;
\item[{\rm ($U_6$)}]
$V_k(x)=U_k(y)$, $V_k(y)=U_k(x)$ and $x\in V_k(y)$ if $y\in V_k(x)$.
\end{itemize}
For $x,y\in \C$ we define the mapping $h_{x,y}\in\aut \C$. Let $z\in \C$. Let's put
\[
h_{x,y}(z)=\begin{cases}
x,& z=y \\
\fuv {V_k(x),V_k(y)}(z),&\text{if $z\in V_k(x)$ for some $k\in\om$}
\end{cases}
\]
It follows from the construction
\begin{itemize}
\item[{\rm ($h_1$)}]
$h_{x,y}\in\aut \C$ for $x,y\in \C$ and $h_{x,y}\in \aut{\C|X}$ for $x,y\in X$;
\item[{\rm ($h_2$)}]
$h_{x,x}=\id_\C$;
\item[{\rm ($h_3$)}]
$h_{x,y} \circ h_{y,z}=h_{x,z}$;
\item[{\rm ($h_4$)}]
$h_{x,y}(y)=x$;
\item[{\rm ($h_5$)}]
$h_{x,y}(U_k(y))=U_k(x)$ and $h_{x,y}(V_k(y))=V_k(x)$.
\item[{\rm ($h_6$)}]
if $x'\in U_k(x)$ and $y'\in U_k(y)$ then
$\restr{h_{x,y}}{V}=\restr{h_{x',y'}}{V}$, where $V=V_k(y)$.
\end{itemize}
Let's put
\[
\pp: \C^3\to \C,\quad (x,y,z) \mapsto h_{x,y}(z).
\]
The condition ($\pp_1$) follows from ($h_2$) and ($h_4$). The condition ($\pp_2$) follows from ($h_3$).
Hence the operation $\pp$ is a strong Mal'tsev operation.

Let us prove that $\pp$ is continuous. Let $x,y,z,u\in \C$, $\pp(x,y,z)=u$, $U$ be a neighborhood of the point $u$. We need to find a neighborhood $W\subset \C^3$ of the point $(z,y,z)$, so that $\pp(W)\subset U$.

Consider the case $y=z$. Then $x=u$.
Then $S=U_k(x)\subset U$ for some $k\in\om$.
Let $V=U_k(y)$ and $W=S\times V \times V$. Let's check $\pp(W)\subset U$. Let $(x',y',z')\in W$. From ($U_5$) and ($h_5$) it follows that $h_{x',y'}(V)=S$. Hence $\pp(x',y',z')\in U$.

Consider the case $y\neq z$. Then $z\in V_m(y)$ for some $m\in \om$. Let $S=U_m(x)$ and $V=U_m(y)$. Since the mapping $h_{x,y}$ is continuous, then $h_{x,y}(Q')\subset U$ for some neighborhood $Q'$ of the point $z$. Let $Q=Q'\cap V_m(y)$ and $W=S\times V \times Q$. Let's check $\pp(W)\subset U$. Let $(x',y',z')\in W$. From ($U_6$) it follows that $\restr{h_{x',y'}}Q=\restr{h_{x,y}}Q$. Hence $\pp(x',y',z')\in U$.

The continuity of $\pp$ is proved. From ($h_1$) it follows that $\pp(X^3)=X$.
\end{proof}

\begin{proposition}\label{prexscsc1}
Let $X$ be a zero-dimensional first countable  space. Then $X^\om\times\C$ is strongly homogeneous.
\end{proposition}
\begin{proof}
Let $Y=X\times\D$. Then the spaces $X^\om\times\C$, $Y^\om$ and $Y^\om \times\D$ are homeomorphic. Therefore, $Y^\om$ has a proper open-closed subspace homeomorphic to $Y^\om$. Proposition 24 (5) \cite{Medini2011} implies that $Y^\om$ is strongly homogeneous.
\end{proof}

\begin{proposition}\label{prexscsc2}
Let $X$ be the \smz/ space. There is a space \smz/ space $Y$ such that $X\times Y$ is strongly homogeneous and strongly rectifiable.
\end{proposition}
\begin{proof}
Let $Y=X^\om\times\C$. Then, by the Proposition \ref{prexscsc1}, $Y$ is strongly homogeneous. The space $X\times Y$ is homeomorphic to $Y$. It remains to apply the Theorem \ref{trexscsc-main}.
\end{proof}

%% file: sec/hsnrtg.tex
\section{Homogeneous spaces that are not right topological groups}\label{sec-hsnrtg}
 
Van Mill in \cite{vanMill1982} constructed a homogeneous space $A$, a subspace of the real numbers, which is not a topological group. In fact, the space $A$ is not a right topological group either.

We call a set $A\subset X$ be a \term{clumsy} if $h(A)\cap A$ is non empty for any homeomorphism $h:X\to X$.

\begin{theorem}\label{thsnrtg1}
Let $X$ be an space and $A\subset X$ be a clumsy subset. If $|A|<|X|$, then $X$ is not a right topological group.
\end{theorem}
\begin{proof}
From the contrary. Take $g\in X \setminus A^{-1}A$. Then $Ag\cap A=\es$. Contradiction.
\end{proof}

Theorem \ref{thsnrtg1} generalizes Lemma 3.4. from \cite{vanMill1982} and the proof of the theorem is much simpler.

There is a countable dense set $E\subset A$ such that $E$ is clumsy in $A$ (Theorem 3.3 \cite{vanMill1982}). Therefore, by the Theorem \ref{thsnrtg1}, the space $A$ from \cite{vanMill1982} is not a right topological group.

In \cite{Banakh2008} Taras Banakh noticed that the zero-dimensional homogeneous first-countable van Douwen compactum \cite{vanDouwen1984} contains a countable clumsy set and, therefore, is not a right topological group.

Let $\T=\P\times \P\cup \Q\times \Q$ is van Douwen space and $\S=\C\times \T$ is van Mill space. Spaces $\T$ and $\S$ are \smz/ strongly homogeneous spaces \cite{vanEngelen1986,vanMill1980}.

\begin{assertion}\label{ahsnrtg1}
$\Q\times \Q$ is clumsy in $\T$.
\end{assertion}
\begin{proof}
Assume the contrary, that is, $h(\Q\times \Q)\cap \Q\times \Q=\es$ for some homeomorphism $h:\T\to \T$.
Let $F=\Q\times \{0\}$. Then $F$ is closed in $\T$ and homeomorphic to $\Q$. Hence $h(F)$ is closed in $\P\times \P$ and homeomorphic to $\Q$. Contradiction.
\end{proof}

The Theorem \ref{thsnrtg1} and Statement \ref{ahsnrtg1} imply

\begin{proposition}\label{phsnrtg1}
$\T$ is not a right topological group.
\end{proposition}

It is not clear about the  van Mill space $\S$ whether the structure a right topological group is allowed on it.

\begin{proposition}[Corollary 5.4 \cite{ArhangelskiiChoban2011}]\label{phsnrtg2}
If a topological group $G$ contains a dense \v Cech complete subspace, then $G$ is \v Cech complete.
\end{proposition}

Recall that Polish spaces are exactly \v Cech complete separable metrizable spaces.

The space $\S$ contains a dense subspace homeomorphic to $\P$ and is not Polish. 
Then $\S^\om$ also contains a dense subspace homeomorphic to $\P$ and is not Polish.
The following sentence follows from the Proposition \ref{phsnrtg2}.

\begin{proposition}\label{phsnrtg3}
$\S^\om$ is not a topological group.
\end{proposition}

%% file: sec/qe.tex
\section{Examples and questions}\label{sec-qe}


\begin{question}(A.V.~Arhangelskii)\label{qqe-1-2}
Let $G$ be a right topological group.

Is it true that if $G$ contains a dense \v Cech complete subspace, then $G$ is \v Cech complete?

Is it true that if $G$ is a separable metrizable and contains a dense Polish subspace, then $G$ is Polish?
\end{question}

\begin{example}\label{eqe1} The van Mill space $\S$ and the van Dawn space $\T$ have the following properties:
\begin{itemize}
\item
are strongly homogeneous \cite{vanMill1980};
\item
absolutely Borel \cite{vanMill1980,vanEngelen1986};
\item
coset sets of some absolutely Borel groups \cite{vanMill2007};
\item
have a dense Polish subspace and are not Polish \cite{vanMill1980};
\item
precompact strongly rectifiable (Theorem \ref{trexscsc-main});
\item
$\T$ is not a right topological group (Proposition \ref{phsnrtg1});
\item
$\S$ is not a topological group (Proposition \ref{phsnrtg2});
\end{itemize}
\end{example}

\begin{question}\label{qqe-1-1}
Is the van Mill space $\S$ a right topological group?
\end{question}

\begin{question}\label{qqe-1}
Let $X\in\sset{\S,\T}$.
\begin{enumerate}
\item
Is the space $X\times \Q$ homeomorphic to a (boolean) (precompact) (right) topological group?
\item
Is the space $X\times \P$ homeomorphic to a (boolean) (precompact) (right) topological group?
\item
Is there a \smz/ space $Y$ such that $X\times Y$ is homeomorphic to a (boolean) (precompact) (right) topological group?
\end{enumerate}
\end{question}

A topological group $G$ is called \term{locally precompact} if there exists a neighborhood $U$ of the neutral element $e$ of $G$ such that $U$ can be covered by ﬁnitely many left and right translates of each neighborhood of $e$ in $G$.
A group is locally precompact iff its completion is locally compact.

\begin{proposition}[{Theorem 8 \cite{Medvedev2013}}]\label{p:med:1} 
Let $G$ be a \smz/ topological group which is not locally precompact. Then the space $G$ is strongly homogeneous.
\end{proposition}

There exist a dense subgroup $G_M$ of the real line \cite[Theorem 5.1]{vanMill1992} and
dense subgroup $G_D$ of the circle \cite{vanDouwen1984}, so $G\in\sset{G_M,G_D}$ satisfies the condition:
\begin{itemize}
\item[]
if $A,B\subset \widehat G$ are open sets and $A\cap G$ is homeomorphic to $B\cap G$, then $\mu(A)=\mu(B)$,
\end{itemize}
where $\widehat G$ is the completion of $G$ (for $G_M$ the completion is the real line, and for $G_D$ the completion is the circle) and $\mu$ is the standard invariant Lebesgue measure on $\widehat G$.
The group $G_D$ is precompact, the group $G_M$ is locally precompact but not precompact.

\begin{example}\label{eqe2}
\smz/ the groups $G_D$ and $G_M$ are not strongly homogeneous, it suffices to take two segments $A,B$ of different lengths, the ends of which do not lie in the group \cite{vanDouwen1984}.
The group $G_M$ is not homeomorphic to a precompact rectifiable space, a finite number of shifts of the set $(-1,1)\cap G_M$ cannot cover $G_M$ (see Lemma \ref{lrexscs1}).
\end{example}

\begin{proposition}[{Corollary 5 \cite{Medvedev2013}}]\label{p:med:1} 
If $X$ is a homogeneous \smz/ space, then $X\times \Q$ is a strongly homogeneous space.
\end{proposition}

\begin{question}\label{qqe1+1}
Let $G\in\sset{G_M,G_D}$.
Which of the spaces listed below are strongly homogeneous:
$G\times \P$, $G\times \C$, $G^\om$?
\end{question}

\begin{question}\label{qqe1+2}
Which of the following spaces are homeomorphic to
(1) precompact group; (2) precompact (strongly) rectifiable space:
$G_M\times \Q$, $G_M\times \P$, $G_M\times \C$, $G_M^\om$?
\end{question}

\begin{question}\label{qqe2}
Let $G\in\sset{G_M,G_D}$.
Which of the spaces listed below are homeomorphic to a Boolean group:
$G$, $G\times \Q$, $G\times P$, $G\times \C$, $G^\om$?
\end{question}

\begin{proposition}[{\cite{Lawrence1998,DowPearl1997}}]\label{p:ldp:1}
Let $X$ be a \smz/ space. Then $X^\om$ is homogeneous.
\end{proposition}

\begin{question}[{Terada, \cite{Terada1993}}]\label{q:terada:1} 
Let $X$ be a \smz/ space.
Is it true that $X^\om$ is a strongly homogeneous space?
\end{question}

\begin{question}\label{qqe3}
Let $X$ be the \smz/ space.
\begin{enumerate}
\item
Which of the following spaces are (1) (strongly) (precompactly) rectifiable spaces; (2) is homeomorphic to a right topological group:
$X^\om$, $X^\om\times \Q$, $X^\om\times \P$, $X^\om\times \C$?
\item
Is the space $X^\om\times \Q$ homeomorphic to a topological (boolean, precompact) group?
\item
Is there a \smz/ space $Y$ such that $X\times Y$ is homeomorphic to a (right) topological group?
\end{enumerate}
\end{question}

\begin{question}\label{qqe4}
Let $X$ be a homogeneous \smz/ space.
\begin{enumerate}
\item
Will the space $X$ be a (strongly) (precompact) rectifiable space?
\item
Is the space $X\times \Q$ homeomorphic to a (precompact) (right) topological group?
\end{enumerate}
\end{question}

\begin{question}\label{qqe5}
In the class \smz/ spaces, distinguish the following classes:
rectifiable; strongly straightened; pre-compact straightening; strongly precompactly rectifiable.
\end{question}